\documentclass[12pt,thmsb]{article}
\usepackage{graphicx}
\usepackage{amsmath}
\usepackage{amsfonts}
\usepackage{amssymb}
\newtheorem{theorem}{Theorem}

\newtheorem{conjecture}[theorem]{Conjecture}
\newtheorem{corollary}[theorem]{Corollary}

\newtheorem{definition}[theorem]{Definition}

\newtheorem{lemma}[theorem]{Lemma}

\newtheorem{proposition}[theorem]{Proposition}
\newtheorem{remark}[theorem]{Remark}

\newenvironment{proof}[1][Proof]{\textbf{#1.} }{\ {\rule{0.5em}{0.5em}}}

\begin{document}

\title{Bijective and general arithmetic codings for Pisot automorphisms
of the torus\footnote{\vskip0.3mm
$\,$ {\em Key words and phrases}: Pisot automorphism,
arithmetic coding, homoclinic point, integral form.}}
\author{Nikita Sidorov\footnote{Supported by the EPSRC grant number
GR/L98923. The author wishes to thank Anatoly Vershik for our fruitful
collaboration in dimension~2 in \cite{SV98} and subsequent stimulating
discussions.}
\\Department of Mathematics, UMIST, P.O. Box 88,
\\Manchester M60 1QD, United Kingdom.
\\E-mail: \texttt{Nikita.A.Sidorov@umist.ac.uk}}
\maketitle
\begin{abstract}
Let $T$ be an algebraic
automorphism of $\mathbb{T}^{m}$ having
the following property: the
characteristic polynomial of its matrix is irreducible over $\mathbb{Q}$, and
a Pisot number $\beta$ is one of its roots. We
define the mapping $\varphi_{\mathbf{t}}$ acting from the two-sided
$\beta$-compactum onto $\mathbb{T}^{m}$ as follows:
\[
\varphi_{\mathbf{t}}(\bar{\varepsilon})=\sum_{k\in\mathbb{Z}}\varepsilon
_{k}T^{-k}\mathbf{t,}%
\]
where $\mathbf{t}$ is a fundamental homoclinic point for $T$, i.e., a point
homoclinic to $\mathbf{0}$ such that the linear span of its orbit is the whole
homoclinic group (provided such a point exists).
We call such a mapping an arithmetic coding of $T$.
This paper is aimed to show that under some natural
hypothesis on $\beta$ (which is apparently satisfied for all Pisot units) the
mapping $\varphi_{\mathbf{t}}$ is bijective a.e. with respect to the Haar
measure on the torus. Besides, we study the case of more general parameters
$\mathbf{t}$, not necessarily fundamental, and relate the number of preimages
of $\varphi_{\mathbf{t}}$ to certain number-theoretic quantities.
We also give several full criteria for $T$ to admit a bijective arithmetic
coding. This work continues the study begun in \cite{SV98} for
the special case $m=2$.
\end{abstract}

\section{Introduction}

Let $T$ be an algebraic automorphism of the torus $\mathbb{T}^{m}$ given by a
matrix $M\in GL(m,\mathbb{Z})$ with the following property: the characteristic
polynomial for $M$ is irreducible over $\mathbb{Q}$, and a Pisot number
$\beta>1$ is its root (we recall that an algebraic integer is called \textit{a
Pisot number }if it is greater than 1 and all its Galois conjugates are less
than 1 in modulus). Since $\det M=\pm1$, $\beta$ is a \textit{unit}, i.e., an
invertible element of the ring $\mathbb{Z[\beta]=Z[\beta}^{-1}\mathbb{]}$. We
will call such an automorphism a \textit{Pisot automorphism}. Note that since
none of the eigenvalues of $M$ lies on the unit circle, $T$ is hyperbolic.

Our goal is to present a symbolic coding of $T$ which, roughly speaking,
reveals not just the structure of $T$ itself but the natural action of the
torus on itself as well. Let us give more precise definitions.

Let $X_{\beta}$ denote the two-sided $\beta$-\textit{compactum}, i.e., the
space of all \textit{admissible} two-sided sequences in the alphabet
$\{0,1,\ldots,[\beta]\}$. More precisely, a representation of an $x\in
\lbrack0,1)$ of the form
\begin{equation}
x=\pi(\varepsilon_{1},\varepsilon_{2},\ldots):=\sum_{1}^{\infty}%
\varepsilon_{k}\beta^{-k} \label{beta-exp}%
\end{equation}
is called the \textit{$\beta$-expansion} of $x$ if the ``digits''
$\{\varepsilon_{k}\}_{1}^{\infty}$ are obtained by means of the greedy
algorithm (similarly to the decimal expanions), i.e., $\varepsilon
_{1}=\varepsilon_{1}(x)=[\beta x],\ \varepsilon_{k}= \varepsilon_{k}%
(x)=[\beta\tau^{k}(x)]$, where $\tau(x)=\{\beta x\}:= \beta
x\operatorname{mod}1$. The set of all possible sequences $\{\{\varepsilon
_{k}(x)\}_{1}^{\infty}:x\in\lbrack0,1)\}$ is called the (one-sided)
\textit{$\beta$-compactum} and denoted by $X_{\beta}^{+}$. A sequence whose
tail is $0^{\infty}$ will be called \textit{finite}.

The $\beta$-compactum can be described more explicitly. Let $1=\sum
_{1}^{\infty}d_{k}^{\prime}\beta^{-k}$ be the expansion of 1 defined as
follows: $d_{1}^{\prime}=[\beta],\ d_{n}^{\prime}=[\beta\tau^{n}1],\;n\geq2$.
If the sequence $\{d_{n}^{\prime}\}$ is not finite, we put $d_{n}\equiv
d_{n}^{\prime}$. Otherwise let $k=\max\,\{j:d_{j}^{\prime}>0\}$, and
$(d_{1},d_{2},\dots):=(\overline{d_{1}^{\prime},\dots,d_{k-1}^{\prime}%
,d_{k}^{\prime}-1})$, where the bar denotes the period of a purely periodic sequence.

We will write $\{x_{n}\}_{1}^{\infty}\prec\{y_{n}\}_{1}^{\infty}$ if
$\{x_{n}\}_{1}^{\infty}\neq\{y_{n}\}_{1}^{\infty}$ and $x_{n}<y_{n}$ for the
smallest $n\geq1$ such that $x_{n}\neq y_{n}$. Then by definition,
\[
X_{\beta}^{+}=\{\{\varepsilon_{n}\}_{1}^{\infty}:(\varepsilon_{n}%
,\varepsilon_{n+1},\dots)\prec(d_{1},d_{2},\dots)\ \text{for all}%
\ n\in\mathbb{N}\}
\]
(see \cite{Pa}). Similarly, we define the \textit{two-sided $\beta$-compactum}
as
\[
X_{\beta}=\{\{\varepsilon_{n}\}_{-\infty}^{\infty}:(\varepsilon_{n}%
,\varepsilon_{n+1},\dots)\prec(d_{1},d_{2},\dots)\ \text{for all}%
\ n\in\mathbb{Z}\}.
\]
Both compacta are naturally endowed with the weak topology, i.e. with the
topology of coordinate-wise convergence, as well as with the natural shifts.
Let the $\beta$-\textit{shift} $\sigma_{\beta}:X_{\beta}\rightarrow X_{\beta}$
act as follows: $\sigma_{\beta}(\bar{\varepsilon})_{k}=\varepsilon_{k+1}$, and
$\sigma_{\beta}^{+}$ be the corresponding one-sided shift on $X_{\beta}^{+}$.
For a Pisot $\beta$ the properties of the $\beta$-shift are well-studied. Its
main property is that it is \textit{sofic}, i.e., is a projection of a
subshift of finite type. In fact this is equivalent to $\{d_{n}\}_{1}^{\infty
}$ being eventually periodic (see, e.g., the review \cite{Bl}).

We extend the $\beta$-expansions to the nonnegative integers in the usual way
(similarly to the decimal expansions).

\begin{lemma}
\label{ep}(see \cite{Ber77}, \cite{Sch80}) Any nonnegative element of the ring
$\mathbb{Z}[\beta]$ has an eventually periodic $\beta$-expansion if $\beta$ is
a Pisot number.
\end{lemma}

There is a natural operation of addition in $X_{\beta}$, namely, if both
sequences $\bar{\varepsilon}$ and $\bar{\varepsilon}^{\prime}$ are
\textit{finite to the left }(i.e., there exists $N\in\mathbb{Z}$ such that
$\varepsilon_{k}=\varepsilon_{k}^{\prime}=0,\;k\leq N$), then by definition,
$\bar{\varepsilon}^{\prime}+\bar{\varepsilon}=\bar{\varepsilon}^{\prime\prime
}$ such that $\sum_{k}\varepsilon_{k}^{\prime\prime}\beta^{-k}=\sum
_{k}(\varepsilon_{k}^{\prime}+\varepsilon_{k})\beta^{-k}$. Later we will show
that under some natural assumption on $\beta$ this operation can be extended
to sequences which are not necessarily finite to the left.

Let $Fin(\beta)$ denote the set of nonnegative $x$'s whose $\beta$-expansions
are finite. Obviously, $Fin(\beta)\subset\mathbb{Z}[\beta]_{+}:=
\mathbb{Z}[\beta]\cap\mathbb{R}_{+}$, but the inverse inclusion does not hold
for an arbitrary Pisot unit.

\begin{definition}
A Pisot unit $\beta$ is called \textbf{finitary} if
\[
Fin(\beta)=\mathbb{Z}[\beta]_{+}.
\]
\end{definition}

A large class of Pisot numbers considered in \cite{FrSo} is known to have this
property. A practical algorithm for checking whether a given Pisot number is
finitary was suggested in \cite{Ak0}. Here is a simple example showing that
not every Pisot unit is finitary. Let $r\geq3$, and $\beta^{2}=r\beta-1$. Then
$X_{\beta}=\{\bar{\varepsilon}:0\leq$ $\varepsilon_{k}\leq r-1,\;(\varepsilon
_{k},\ldots,\varepsilon_{k+n})\neq(r-1,r-2,\ldots,r-2,r-1),\;k\in
\mathbb{Z},n\geq1\}$ and $1-\beta^{-1}=(r-2)\beta^{-1}+(r-2)\beta^{-2}+\ldots
$, i.e., $\beta$ is not finitary.

\begin{definition}
A Pisot unit $\beta$ is called \textbf{weakly finitary}\textit{ if for any
}$\delta>0$ and any $x\in\mathbb{Z}[\beta]_{+}$ there exists $f\in
Fin(\beta)\cap\lbrack0,\delta)$ such that $x+f\in Fin(\beta)$ as well. \label{wf}
\end{definition}

This condition was considered in the recent work by Sh.~Akiyama \cite{Ak}, in
which the author shows that the boundary of the natural sofic tiling generated
by a weakly finitary Pisot $\beta$ has Lebesgue measure zero (moreover, these
conditions are actually equivalent). The author is grateful to Sh.~Akiyama for
drawing his attention to this paper and for helpful discussions.

A slightly weaker (but possibly equivalent) condition
\[
\mathbb{Z}[\beta]=Fin(\beta)-Fin(\beta)
\]
together with the finiteness of $\{d_{n}^{\prime}\}$ was used in the recent
Ph.D. dissertation \cite{Hol} to show that the spectrum of the Pisot
substitutional dynamical system%
\[
0\rightarrow10^{d_{1}^{\prime}},1\rightarrow20^{d_{2}^{\prime}},\ldots
,l-2\rightarrow(l-1)0^{d_{l-1}^{\prime}},l-1\rightarrow0^{d_{l}^{\prime}},
\]
(where $l=\max\;\{n:d_{n}^{\prime}\neq0\}$), is purely discrete. This claim is
a generalization of the corresponding result for a finitary $\beta$ from
\cite{BSolomyak} (see also \cite{MSolomyak}).

\begin{conjecture}
\label{conject}Any Pisot unit is weakly finitary.
\end{conjecture}

To support this conjecture, we are going to explain how to verify that a
particular Pisot unit is weakly finitary. Firstly, one needs to describe all
the elements of the set
\begin{equation}
Z_{\beta}=\{\alpha\in\mathbb{Z}[\beta]\cap\lbrack0,1):\alpha\;\text{has a
purely periodic }\beta\text{-expansion}\}. \label{Zbeta}%
\end{equation}

\begin{lemma}
(see \cite{Ak}). The set $Z_{\beta}$ is finite.
\end{lemma}

\begin{proof}
The sketch of the proof is as follows: basically, the claim follows from
Lemma~\ref{structure}, which implies that the denominator of any $\alpha
\in\mathcal{P}_{\beta}$ in the standard basis of
$\mathbb{Q}(\beta)$ is uniformly bounded, whence the period of the $\beta
$-expansion of $\alpha$ is bounded as well.
\end{proof}

Therefore, we have a finite collection of numbers $\{\sum_{j=0}^{m-1}%
y_{j}\beta^{j}:|y_{j}|\leq q\}$ to ``check for periods'' (here $q$ is the
denominator of $\xi_{0}$ defined by (\ref{P_beta}) in the standard basis of
the ring). Next, it is easy to see that if suffices to check that Definition
\ref{wf} holds for any $x=\alpha\in Z_{\beta}$ (see \cite{Ak}). Moreover, we
can confine ourselves to the case $f\in Fin(\beta)\cap\lbrack\beta^{-2p}%
,\beta^{-p})$, where $p$ is the period of $\alpha$. Indeed, if such an $f$
exists, $\beta^{-p}f$ will do as well, and we will be able to make $f$
arbitrarily small. All known examples of Pisot units prove to be weakly finitary.

We will need the following technical result.

\begin{lemma}
\label{xi}A Pisot unit is weakly finitary if and only if the following
condition is satisfied: there exists $\eta=\eta(\beta)\in(0,1)$ such that
\textit{for any }$\delta>0$ and any $x\in\mathbb{Z}[\beta]_{+}$ there exists
$f\in Fin(\beta)\cap\lbrack\eta\delta,\delta)$ such that $x+f\in Fin(\beta)$
as well.
\end{lemma}

\begin{proof}
It suffices to show that if $\beta$ is weakly finitary, then $\eta$ in
question does exist. Let $\beta$ be weakly finitary; then for any $\alpha\in
Z_{\beta}$ there exists $f_{\alpha}\in Fin(\beta)$ such that $\alpha
+f_{\alpha}\in Fin(\beta)$. Let $\alpha$ has the $\beta$-expansion
$(\overline{\alpha_{1},\ldots,\alpha_{p}})$ and $\alpha^{\ast}=\sum_{1}%
^{p}\alpha_{j}\beta^{-j}.$ Without loss of generality we may regard $p$ to be
greater than the preperiod $+$ the period of the sequence $\{d_{n}%
\}_{1}^{\infty}$ (as $p$ is not necessarily the \textit{smallest} period of
$\alpha$). Since $f_{\alpha}$ can be made arbitrarily small, we may fix it
such that%
\begin{equation}
\alpha+f_{\alpha}<\alpha^{\ast}+\beta^{-p}\alpha. \label{alpha+f}%
\end{equation}
Put $\eta:=\min\;\{f_{\alpha}:\alpha\in Z_{\beta}\}$.

Let $x\in\mathbb{Z}[\beta]_{+}$. By Lemma~\ref{ep} the $\beta$-expansion of
$x$ is eventually periodic, and splitting it into the preperiodic and periodic
parts, we have $x=x_{0}+\beta^{-k}\alpha,\;x_{0}\in Fin(\beta),k\in
\mathbb{Z},\alpha\in Z_{\beta}$. Let for simplicity of notation $k=0$ (the
whole picture is shift-invariant). It will suffice to check the condition for
$\delta=\delta_{n}=\beta^{-pn}$. Put $f=f_{n}:=\beta^{-pn}f_{\alpha}$. Then
\begin{equation}
x+f=(x_{0}+\alpha^{\ast}+\beta^{-p}\alpha^{\ast}+\cdots+\beta^{-(n-1)p}%
\alpha^{\ast})+\beta^{-pn}(\alpha+f_{\alpha}) \label{x+f}%
\end{equation}
The expression in brackets in (\ref{x+f}) belongs to $Fin(\beta)$ and so does
the second term. In view of (\ref{alpha+f}) and the definition of $X_{\beta}$
the whole sum in (\ref{x+f}) belongs to $Fin(\beta)$ as well, because by our
choice of $p$ we have necessarily $(\alpha_{1},\ldots,\alpha_{p})\prec
(d_{1},\ldots,d_{p})$. Since $Z_{\beta}$ is finite and the construction
depends on $\alpha$ only, we are done.
\end{proof}

\section{Formulation of the main result and first steps of the proof}

We recall that the hyperbolicity of $T$ implies that it has the stable and
unstable foliation and consequently the set of homoclinic points. More
precisely, a point $\mathbf{t\in}\mathbb{T}^{m}$ is called \textit{homoclinic
to zero }or simply \textit{homoclinic }if $T^{n}\mathbf{t}\rightarrow
\mathbf{0}$ as $n\rightarrow\pm\infty$ (as is well known, the convergence to
$\mathbf{0}$ in this case will be at exponential rate). In other terms, a
homoclinic point $\mathbf{t}$ must belong to the intersection of the leaves of
the stable foliation $L_{s}$ and the unstable foliation $L_{u}$ passing
through $\mathbf{0}$. Let $\mathcal{H}(T)$ denote the set of all homoclinic
points for $T$; obviously, $\mathcal{H}(T)$ is a group under addition. In
\cite{Ver} it was shown that every homoclinic point can be obtained by
applying the following procedure: take a point $\mathbf{n}\in\mathbb{Z}^{m}$
and project it onto $L_{u}$ along $L_{s}$. Let $\mathbf{s}$ denote this
projection; finally, project $\mathbf{s}$ onto the torus by taking the
fractional parts of all its coordinates. The correspondence $\mathbf{n}%
\leftrightarrow\mathbf{s}\leftrightarrow\mathbf{t}$ is one-to-one. We will
call $\mathbf{s=s(t)}$ the $\mathbb{R}^{m}$-\textit{coordinate }of a
homoclinic point $\mathbf{t}$ and $\mathbf{n}$ the $\mathbb{Z}^{m}%
$-\textit{coordinate} of $\mathbf{t}$. Note that since $T$ is a Pisot
automorphism, we have $\dim L_{u}=1,\;\dim L_{s}=m-1$.

We wish to find an \textit{arithmetic coding} $\varphi$ of $T$ in the
following sense: we choose $X_{\beta}$ as a symbolic compact space and impose
the following restrictions on a map $\varphi:X_{\beta}\rightarrow
\mathbb{T}^{m}$:

\begin{enumerate}
\item $\varphi$ is continuous and bounded-to-one;

\item $\varphi\sigma_{\beta}=T\varphi$;

\item $\varphi(\bar{\varepsilon}+\bar{\varepsilon}^{\prime})=\varphi
(\bar{\varepsilon})+\varphi(\bar{\varepsilon}^{\prime})$ for any pair of
sequences finite to the left.
\end{enumerate}

In \cite{SV98} it was shown that if $m=2$, then there exists $\mathbf{t\in
}\mathcal{H}(T)$ such that $\varphi=\varphi_{\mathbf{t}}:X_{\beta}%
\rightarrow\mathbb{T}^{m}$:%
\begin{equation}
\varphi_{\mathbf{t}}(\bar{\varepsilon})=\sum_{k\in\mathbf{Z}}\varepsilon
_{k}T^{-k}\mathbf{t}. \label{phi}%
\end{equation}
The proof for an arbitrary $m$ is basically the same, and we will omit it. Our
primary goal is to find an arithmetic coding that is bijective a.e. Let us
make some remarks.

Note that the idea of using homoclinic points to ``encode'' ergodic toral
automorphisms had been suggested by A.~Vershik in \cite{Mittag} for $M=%
\begin{pmatrix}
1 & 1\\
1 & 0
\end{pmatrix}
$ and was later developed for a more general context in numerous works -- see
\cite{Ver}, \cite{Leborgne}, \cite{SV97}, \cite{SV98}, \cite{Sch}. The choice
of $X_{\beta}$ as a ``coding space'' is special in the case in question;
indeed, the topological entropy of the shift $\sigma_{\beta}$ is known to be
$\log\beta$ and so is the entropy of $T$. In a more general context (for
example, if $M$ has two eigenvalues outside the unit disc) it is still
unclear, which compactum might replace $X_{\beta}$. Indeed, since $\varphi$ is
bounded-to-one, the topological entropy of the subshift on this compactum must
have the same topological entropy as $T$, i.e., $\log\prod_{|\beta_{j}%
|>1}|\beta_{j}|$, where $\beta_{j},\ j=1,\dots,m,$ are the conjugates of
$\beta$, and there is apparently no natural subshift associated with $\beta$
which has this entropy. However, it is worth noting that the
\textit{existence} of such compacta in different settings has been shown in
\cite{Ver}, \cite{KenVer}, \cite{Sch}.

Return to our context. The mapping $\varphi_{\mathbf{t}}$ defined by
(\ref{phi}) is indeed well defined and continuous, as the series (\ref{phi})
converges at exponential rate. Furthermore, since $T^{k}\mathbf{t=\beta
t}\operatorname{mod}\mathbb{Z}^{m}$, we have by continuity $\varphi
_{\mathbf{t}}\sigma_{\beta}=T\varphi_{\mathbf{t}}$, i.e., $\varphi
_{\mathbf{t}}$ does semiconjugate the shift and a given automorphism $T$.

We will call $\varphi_{\mathbf{t}}$ a \textit{general arithmetic coding} for
$T$ (parametrised by a homoclinic point $\mathbf{t}$).

\begin{lemma}
For any choice of $\mathbf{t}$ the mapping $\varphi_{\mathbf{t}}$ is bounded-to-one.\label{bounded}
\end{lemma}

\begin{proof}
Let $\left\|  \cdot\right\|  $ denote the distance to the closest integer,
$\mathbf{s}$ be the $\mathbb{R}^{m}$-coordinate of $\mathbf{t}$ and $\tilde
{T}$ denote the linear transformation of $\mathbb{R}^{m}$ defined by the
matrix $M$. Let $\varphi_{N,\mathbf{t}}$ be the mapping acting from $X_{\beta
}$ into $\mathbb{R}^{m}$ by the formula
\[
\phi_{N,\mathbf{t}}(\bar{\varepsilon}):=\sum_{-N}^{N}\varepsilon_{k}(\tilde
{T})^{-k}\mathbf{s}.
\]
Then by (\ref{phi}),
\[
\varphi_{\mathbf{t}}(\bar{\varepsilon})=\lim_{N\rightarrow+\infty}
(\phi_{N,\mathbf{t}}(\bar{\varepsilon})\operatorname{mod}\;\mathbb{Z}^{m}),
\]
where $(x_{1},\dots,x_{m})\operatorname{mod}\;\mathbb{Z}^{m}= (\{x_{1}%
\},\dots,\{x_{m}\})$. Therefore, it suffices to show that the diameters of the
sets $\phi_{N,\mathbf{t}}(X_{\beta})$ are uniformly bounded for all $N$. We
have (recall that $0\leq\varepsilon_{k}\leq\lbrack\beta]$):
\begin{align*}
\max\,\{\|\phi_{N,\mathbf{t}}(\overline\varepsilon)\|:\overline\varepsilon
\in X_{\beta}\}  &  \leq\lbrack\beta]\left\|
\sum_{-N}^{N}(\tilde{T})^{-k}\mathbf{s}\right\|  \leq\lbrack\beta]\sum
_{-N}^{N}\left\|  (\tilde{T})^{-k}\mathbf{s}\right\| \\
&  \leq\text{const}\cdot\sum_{0}^{N}\theta^{k}<\infty,
\end{align*}
where $\theta\in(0,1)$ is the maximum of the absolute values of the conjugates
of $\beta$ that do not coincide with $\beta$. This proves the lemma.
\end{proof}

Let the characteristic equation for $\beta$ be
\[
\beta^{m}=k_{1}\beta^{m-1}+k_{2}\beta^{m-2}+\cdots+k_{m}%
\]
and $T_{\beta}$ denote the toral automorphism given by the \textit{companion
matrix }$M_{\beta}$ for $\beta$, i.e.,%
\[
M_{\beta}=\left(
\begin{array}
[c]{ccccc}%
k_{1} & k_{2} & \ldots &  k_{m-1} & k_{m}\\
1 & 0 & \ldots & 0 & 0\\
0 & 1 & \ldots & 0 & 0\\
\ldots & \ldots & \ldots & \ldots & \ldots\\
0 & 0 & \ldots & 1 & 0
\end{array}
\right)  .
\]
We first assume the following conditions to be satisfied:

\begin{enumerate}
\item $T$ is algebraically conjugate to $T_{\beta}$, i.e., there exists a
matrix $C\in GL(m,\mathbb{Z})$ such that $CM=M_{\beta}C$ (notation: $T\sim
T_{\beta}$).

\item  A homoclinic point $\mathbf{t}$ is \textit{fundamental}, i.e.,
$\left\langle T^{n}\mathbf{t}\mid n\in\mathbb{Z}\right\rangle =\mathcal{H}(T)$.

\item $\beta$ is weakly finitary.
\end{enumerate}

The notion of fundamental homoclinic point for general actions of expansive
group automorphisms was introduced in \cite{LindSch} (see also \cite{Sch}).

\begin{remark}
Note that the second condition implies the first, as the mere existence of a
fundamental homoclinic point means that $T\sim T_{\beta}$ (see
Theorem~\ref{bac} below). Conversely, if $T\sim T_{\beta}$, then there is
always a fundamental homoclinic point for $T$. Indeed, let $\mathbf{n}%
_{0}=(0,0,\ldots,0,1)$ be the $\mathbb{Z}^{m}$-coordinate of $\mathbf{t}_{0}$.
Then $\mathbf{t}_{0}$ is a fundamental for $T_{\beta}$ and if $CM=M_{\beta}C$,
then $C^{-1}\mathbf{t}_{0}$ is fundamental for $T$. \label{conj}
\end{remark}

Now we are ready to formulate the main theorem of the present paper.

\begin{theorem}
Provided the above conditions are satisfied, the mapping $\varphi_{\mathbf{t}%
}$ defined by (\ref{phi}) is bijective a.e. with respect to the Haar measure
on the torus. \label{main}
\end{theorem}

\smallskip\noindent\textbf{Remark.} In \cite{SV98} the claim of the theorem
was shown for $m=2$. We wish to follow the line of exposition of that paper,
though it is worth stressing that our approach will be completely different
(rather arithmetic than geometric). In \cite{Sch} this claim was proven for
any finitary $\beta$ and it was conjectured that it holds for any Pisot
automorphism satisfying conditions 1 and 2 above. We give further support for
this conjecture, as Theorem~\ref{main} implies that we actually reduced it to
a general number-theoretic conjecture verifiable for any given Pisot unit
$\beta$ (see Conjecture~\ref{conject}).

The rest of the section as well as the next section will be devoted to the
proof of Theorem~\ref{main}; in the last section we will discuss the case when
conditions~1 and 2 are not necessarily satisfied.

We are going to need the following number-theoretic claim. Let
\[
\mathcal{P}_{\beta}:=\{\xi:\left\|  \xi\beta^{n}\right\|  \rightarrow
0,\;n\rightarrow+\infty\}.
\]
It is obvious that $\mathcal{P}_{\beta}$ is a group under addition.

\begin{lemma}
\label{structure}There exists $\xi_{0}\in\mathbb{Q}(\beta)\setminus
\mathbb{Z}[\beta]$ such that
\begin{equation}
\mathcal{P}_{\beta}=\xi_{0}\cdot\mathbb{Z}[\beta].\label{P_beta}%
\end{equation}
\end{lemma}

\begin{proof}
By the well-known result, for any Pisot $\beta,\ \xi\in\mathcal{P}_{\beta
}\Leftrightarrow\xi\in\mathbb{Q}(\beta)$ and $Tr(\beta^{k}\xi)\in
\mathbb{Z},\;k\geq k_{0}$ (where $Tr(\varsigma)$ denotes the trace of an
element $\varsigma$ of the extension $\mathbb{Q}(\beta)$, i.e., the sum of all
its Galois conjugates) -- see, e.g., \cite{Cas}. Since $\beta$ is a unit,
$Tr(\varsigma)\in\mathbb{Z}$ implies $Tr(\beta^{-1}\varsigma)\in\mathbb{Z}$,
whence
\begin{equation}
\mathcal{P}_{\beta}:=\{\xi\in\mathbb{Q}(\beta):Tr(a\xi)\in\mathbb{Z\;\forall
}a\in\mathbb{Z}[\beta]\}. \label{dual}%
\end{equation}
Thus, if we regard $\mathbb{Z}[\beta]$ as a lattice over $\mathbb{Z}$, then by
(\ref{dual}), $\mathcal{P}_{\beta}$ is by definition the dual lattice for
$\mathbb{Z}[\beta]$. Hence by the well known ramification theorem (see, e.g.,
\cite[Chapter~III]{FrolTa}) the equality~(\ref{P_beta}) follows with $\xi
_{0}=1/g^{\prime}(\beta)$, where $g(x)=x^{m}-k_{1}x^{m-1}-\cdots-k_{m}$.
\end{proof}

We will divide the proof of the main theorem into several steps.

\noindent\textbf{Step 1 (description of the homoclinic group).}

\begin{lemma}
Any homoclinic point $\mathbf{t}$ for $T_{\beta}$ has the $\mathbb{R}^{m}%
$-coordinate%
\begin{equation}
\mathbf{s(t)}=\xi_{0}u(1,\beta^{-1},\ldots,\beta^{-m+1}), \label{u}%
\end{equation}
where $u\in\mathbb{Z}[\beta]$.
\end{lemma}

\begin{proof}
We have $M_{\beta}\bar{v}_{\beta}=\beta\bar{v}_{\beta}$, where $\bar{v}%
_{\beta}=(1,\beta^{-1},\ldots,\beta^{-m+1})$. As was mentioned above, the
dimension of the unstable foliation $L_{u}$ is 1, whence $\mathbf{s}%
(\mathbf{t})=k\bar{v}_{\beta}$, and since $T_{\beta}^{n}\mathbf{t\rightarrow
0}$, we have $\left\|  k\beta^{n}\right\|  \rightarrow0$, i.e., $k\in
\mathcal{P}_{\beta}$. Now the claim of the lemma follows from (\ref{P_beta}).
\end{proof}

Let $\mathcal{U}_{\beta}$ denote the group of units ($=$ invertible elements)
of the ring $\mathbb{Z}[\beta]$.

\begin{lemma}
\label{11} There is a one-to-one correspondence between the group
$\mathcal{U}_{\beta}$ and the set of fundamental homoclinic points for
$T_{\beta}$. Namely, if $\mathbf{t}$ is fundamental, then $u$ in (\ref{u}) is
a unit and vice versa.
\end{lemma}

\begin{proof}
Suppose $\mathbf{t}$ is fundamental. Then the homoclinic point $\mathbf{t}%
_{0}$ whose $\mathbb{R}^{m}$-coordinate is $\mathbf{s}_{0} =(\xi_{0},\xi
_{0}\beta^{-1},\ldots,\xi_{0}\beta^{-m+1})$ can be represented as a finite
linear integral combination of the powers $T^{k}\mathbf{t}$, i.e.,
\[
\xi_{0}(1,\beta^{-1},\ldots,\beta^{-m+1})=\sum_{k}e_{k}\beta^{k}\xi
_{0}u(1,\beta^{-1},\ldots,\beta^{-m+1}),
\]
whence $u\sum_{k}e_{k}\beta^{k}=1$. Therefore, $u$ is invertible in the ring
$\mathbb{Z}[\beta]$.

Conversely, if $u\in\mathcal{U}_{\beta}$, then using the same method, we show
that the claim of the lemma follows from the fact that the equation
$ux=u^{\prime}$ always has the solution in $\mathbb{Z}[\beta]$, namely,
$x=u^{-1}u^{\prime}$.
\end{proof}

\noindent\textbf{Step 2 (reduction to} $T=T_{\beta}$\textbf{).} To prove
Theorem \ref{main}, we may without loss of generality assume $T=T_{\beta}$.
Indeed, suppose $M=C^{-1}M_{\beta}C$, where $C\in GL(m,\mathbb{Z})$. Then
there is a natural one-to-one correspondence between $\mathcal{H}(T)$ and
$\mathcal{H}(T_{\beta})$, namely, $\mathbf{t\in}\mathcal{H}(T)\Leftrightarrow
C\mathbf{t\in}\mathcal{H}(T_{\beta})$. Furthermore, if $\varphi_{\mathbf{t}}$
is bijective a.e., then so is $\varphi_{C\mathbf{t}}$, as $\varphi
_{C\mathbf{t}}=C\varphi_{\mathbf{t}}$.

So, we assume first that $T=T_{\beta}$, and $\mathbf{t}$ is a general
fundamental homoclinic point for $T_{\beta}$ given by (\ref{u}). In this case
the formula (\ref{phi}) becomes%
\[
\varphi_{\mathbf{t}}(\bar{\varepsilon})=\lim_{N\rightarrow+\infty}\sum
_{k=-N}^{+\infty}\varepsilon_{k}\beta^{-k}\left(
\begin{array}
[c]{c}%
\xi_{0}u\\
\xi_{0}u\beta^{-1}\\
\vdots\\
\xi_{0}u\beta^{-m+1}%
\end{array}
\right)  \operatorname{mod}\;\mathbb{Z}^{m}.\label{phi2}%
\]

\noindent\textbf{Step 3 (the preimage of \textbf{{0})}}. Let $Z_{\beta}$ be
defined by (\ref{Zbeta}).

\begin{lemma}
The preimage of $\mathbf{0}$ can be described as follows:%
\[
\mathcal{O}_{\beta}:=\varphi_{\mathbf{t}}^{-1}(\mathbf{0})=\{\bar{\varepsilon
}\in X_{\beta}:\bar{\varepsilon}\text{ is purely periodic and }\sum
_{1}^{\infty}\varepsilon_{j}\beta^{-j}\in Z_{\beta}\}.
\]
\end{lemma}

\begin{proof}
By Lemma \ref{bounded}, $\mathcal{O}_{\beta}$ is finite and since it is
shift-invariant, it must contain purely periodic sequences only. Let
$\alpha=\sum_{1}^{\infty}\varepsilon_{j}\beta^{-j}$. Then by (\ref{phi2}),
$\left\|  \alpha u\xi_{0}\beta^{n}\right\|  \rightarrow0$ as $n\rightarrow
\infty$, whence from (\ref{P_beta}), therefore, $\alpha u\in\mathbb{Z[\beta]}%
$, and $\alpha\in\mathbb{Z[\beta]}$, because $u\in\mathcal{U}_{\beta}$.
\end{proof}

\noindent\textbf{Step 4 (description of the full preimage of any point of the
torus). }We are going to show that $\varphi_{\mathbf{t}}$ is ``linear'' in the
sense that for any two sequences $\bar{\varepsilon},\bar{\varepsilon}^{\prime
}\in\varphi_{\mathbf{t}}^{-1}(x)$ their ``difference'' will belong to
$\mathcal{O}_{\beta}$. More precisely, let $\varepsilon^{(N)}$ denote both the
sequence $(\ldots,0,0,\ldots,0,\varepsilon_{-N},\varepsilon_{-N+1},\ldots)$
and its ``value'' $e^{(N)}:=\sum_{k=-N}^{+\infty}\varepsilon_{k}\beta^{-k}$

\begin{lemma}
If $\varphi_{\mathbf{t}}(\bar{\varepsilon})=\varphi_{\mathbf{t}}%
(\bar{\varepsilon}^{\prime})$, then for any $N\geq1$ there exists $\alpha\in
Z_{\beta}$ such that
\[
|e^{(N)}-(e^{\prime})^{(N)}|=\beta^{N}\alpha.
\]
\end{lemma}

\begin{proof}
Let $\mathcal{E}$ denote the set of all partial limits of the collection of
sequences $|\varepsilon^{(N)}-(\varepsilon^{\prime})^{(N)}|,N\geq1$, where
$|\varepsilon^{(N)}-(\varepsilon^{\prime})^{(N)}|$ is the sequence
$(\ldots,0,0,\ldots,0,\allowbreak\varepsilon_{-N}^{\prime\prime}%
,\allowbreak\varepsilon_{-N+1}^{\prime\prime},\ldots)$ whose ``value'' is
$|e^{(N)}-(e^{\prime})^{(N)}|$. It suffices to show that $\mathcal{E}%
\subset\mathcal{O}_{\beta}$. Let $\bar{\delta}\in\mathcal{E}$; by definition,
there exists a sequence of positive integers $\{N_{k}\}$ such that
$\delta^{(N_{k})}=\left|  (\varepsilon^{\prime})^{(N_{k})}-\varepsilon
^{(N_{k})}\right|  ,\;k=1,2,\ldots$ Then $\varphi_{\mathbf{t}}(\bar{\delta
})=\lim_{k\rightarrow\infty}\varphi_{\mathbf{t}}(\delta^{(N_{k})})=\mathbf{0}%
$, and we are done.
\end{proof}

Therefore, if $\bar{\varepsilon}\in\varphi_{\mathbf{t}}^{-1}(x)$ for some
$x\in\mathbb{T}^{m}$, then we know that to obtain any $\bar{\varepsilon
}^{\prime}\in\varphi_{\mathbf{t}}^{-1}(x)$, one may take one of the partial
limits of the sequence $\{\varepsilon^{(N)}+\beta^{N}\alpha\}$ for $\alpha\in
Z_{\beta},$ perhaps, depending on $N$. We will write
\begin{equation}
\bar{\varepsilon}\sim\bar{\varepsilon}^{\prime}\mbox{ iff }
\varphi_{\mathbf{t}}(\bar{\varepsilon})=
\varphi_{\mathbf{t}}(\bar{\varepsilon}^{\prime}).
\label{er}
\end{equation}

\smallskip\noindent\textbf{Conclusion.} Thus, we reduced the proof of Theorem
\ref{main} to a certain claim about the two-sided $\beta$-compactum. \smallskip

Basically, our goal now is to show that the procedure described above will not
change an arbitrarily long tail of a generic sequence $\bar{\varepsilon}\in
X_{\beta}$ and therefore, will not change $\bar{\varepsilon}$ itself.

\section{Final steps of the proof and examples}

Let $\mu_{\beta}$ denote the measure of maximal entropy for the shift
$(X_{\beta},\sigma_{\beta})$, and $\mu_{\beta}^{+}$ be its one-sided analog.
We wish to prove that%
\begin{equation}
\mu_{\beta}\{\bar{\varepsilon}\in X_{\beta}:\#[\bar{\varepsilon}]=1\}=1,
\label{mu=1}%
\end{equation}
where $[\bar{\varepsilon}]=\{\bar{\varepsilon}^{\prime}\in X_{\beta}%
:\bar{\varepsilon}^{\prime}\sim\bar{\varepsilon}\}$.

\noindent\textbf{Step 5 (estimation of the measure of the ``bad'' set).} We
will need some basic facts about the measure $\mu_{\beta}$. For technical
reasons we prefer to deal with its one-sided analog $\mu_{\beta}^{+}$.

\begin{lemma}
\label{below} There exists a constant $C_{1}=C_{1}(\beta)\in(0,1)$ such that
for any $n\geq2$ and any $(i_{1},i_{2},\ldots)\in X_{\beta}^{+}$,%
\[
\mu_{\beta}^{+}\;(\varepsilon_{n}=i_{n}\mid\varepsilon_{n-1}=i_{n-1}%
,\ldots,\varepsilon_{1}=i_{1})\geq C_{1}.
\]
\end{lemma}

\begin{proof}
Let the mapping $\pi:X_{\beta}^{+}\rightarrow\lbrack0,1)$ be given by
formula~(\ref{beta-exp}) and $m_{\beta}^{+}=\pi(\mu_{\beta}^{+})$. Let
$C_{n}(\bar{\varepsilon})=(\varepsilon_{n}=i_{n},\varepsilon_{n-1}%
=i_{n-1},\ldots,\varepsilon_{1}=i_{1})\subset X_{\beta}^{+}$ and $\Delta
_{n}(\bar{\varepsilon})=\pi(C_{n}(\bar{\varepsilon})).$ The Garsia Separation
Lemma \cite{Garsia} says that there exists a constant $K=K(\beta)>0$ such that
if $\bar{\varepsilon}$ and $\bar{\varepsilon}^{\prime}$ are two sequences in
$X_{\beta}^{+}$ and $\sum_{k=1}^{n}\varepsilon_{k}\beta^{-k}\neq\sum_{k=1}%
^{n}\varepsilon_{k}^{\prime}\beta^{-k}$, then $\left|  \sum_{k=1}%
^{n}(\varepsilon_{k}-\varepsilon_{k}^{\prime})\beta^{-k}\right|  \geq
K\beta^{-n}$. Hence
\[
K\leq\beta^{n}\mathcal{L}_{1}(\Delta_{n}(\bar{\varepsilon}))\leq1,
\]
where $\mathcal{L}_{1}$ denotes the Lebesgue measure on $[0,1]$. Since for any
$\beta>1,$ $m_{\beta}^{+}$ is equivalent to $\mathcal{L}_{1}$ and the
corresponding density is uniformly bounded away from $0$ and $\infty$ (see
\cite{Re}), we have for some $K^{\prime}>1$,
\[
1/K^{\prime}\leq\beta^{n}m_{\beta}^{+}(\Delta_{n}(\bar{\varepsilon}))\leq
K^{\prime},
\]
whence by the fact that $\pi$ is one-to-one except for a countable set of
points,
\[
1/K^{\prime}\leq\beta^{n}\mu_{\beta}^{+}(C_{n}(\bar{\varepsilon}))\leq
K^{\prime}%
\]
and the claim of the lemma holds with $C_{1}=(\beta K^{\prime})^{-2}$.
\end{proof}

There is a natural arithmetic structure on $X_{\beta}^{+}$: the sum of two
sequences $\bar{\varepsilon}$ and $\bar{\varepsilon}^{\prime}$ is defined as
the sequence equal to the $\beta$-expansion of the sum $\{\sum_{1}^{\infty
}(\varepsilon_{k}+\varepsilon_{k}^{\prime})\beta^{-k}\}$. Let $X_{\beta}%
^{(n)}$ denote the set of finite words of length $n$ that are extendable to a
sequence in $X_{\beta}^{+}$ by writing noughts at all places starting with
$n+1$. We will sometimes identify $X_{\beta}^{(n)}$ with the set
$Fin_{n}(\beta):=\{\bar{\varepsilon}:\varepsilon_{k}\equiv0,\ k\geq n+1\}$.

By the sum $(\varepsilon_{1},\varepsilon_{2},\ldots,\varepsilon_{n}%
)+\bar{\varepsilon}^{\prime}$, we will imply $(\varepsilon_{1},\varepsilon
_{2},\ldots,\varepsilon_{n},0,0,\ldots)+\bar{\varepsilon}^{\prime}$. In
\cite{FrSo} it was shown that there exists a natural $L_{1}=L_{1}(\beta)$ such
that if $\bar{\varepsilon}\in Fin_{n}(\beta),\;\bar{\varepsilon}^{\prime}\in
Fin_{n}(\beta)$ and $\bar{\varepsilon}+\bar{\varepsilon}^{\prime}\in
Fin(\beta)$, then $\bar{\varepsilon}+\bar{\varepsilon}^{\prime}\in
Fin_{n+L_{1}}(\beta)$.

Recall that by Lemma~\ref{xi} there exists $\eta=\eta(\beta)\in(0,1)$ such
that the quantity $f$ in Definition~\ref{wf} can be chosen in $(\eta
\delta,\delta)$ instead of $(0,\delta)$. We set
\[
L_{2}:=\frac{\log(1/\eta)}{\log\beta}.
\]
Let $L:=\max\;\{L_{1},L_{2}\}$. We can reformulate the hypothesis that $\beta$
is weakly finitary as follows ($\bar{\alpha}$ denotes the $\beta$-expansion of
$\alpha$):%
\begin{align}
&\mbox{for any}\,(\varepsilon_{1},\varepsilon_{2},\ldots,\varepsilon_{n})\in
X_{\beta}^{(n)}\,\mbox{there exists}\,
(\varepsilon_{n+1},\ldots,\varepsilon_{n+L})\in
X_{\beta}^{(L)}:\label{WF}\\
&  (\varepsilon_{1},\ldots,\varepsilon_{n+L})\in X_{\beta}^{(n+L)},
\; \bar{\alpha}+(\varepsilon_{1},\ldots,
\varepsilon_{n+L})\in Fin(\beta)\,\mbox{for all}\ \alpha\in Z_{\beta}.\nonumber
\end{align}
A direct consequence of Lemma~\ref{below} is

\begin{corollary}
For any $(i_{1},i_{2},\ldots)\in X_{\beta}^{+}$,
\begin{equation}
\frac{\mu_{\beta}^{+}(\varepsilon_{n+L}=i_{n+L},\varepsilon_{n+L-1}%
=i_{n+L-1},\ldots,\varepsilon_{1}=i_{1})}{\mu_{\beta}^{+}(\varepsilon
_{n}=i_{n},\varepsilon_{n-1}=i_{n-1},\ldots,\varepsilon_{1}=i_{1})}\geq
C_{2}=C_{1}^{L}. \label{C_1}%
\end{equation}
\end{corollary}

\begin{lemma}
\label{glue} If $(i_{1},\ldots,i_{n})\in X_{\beta}^{(n)}$ and $(j_{1}%
,\ldots,j_{k})\in X_{\beta}^{(k)}$, then $(i_{1},\ldots,i_{n}%
,0,0,0,0,\allowbreak j_{1}, \ldots,j_{k})\in X_{\beta}^{(n+k+4)}$.
\end{lemma}

\begin{proof}
The claim follows from the definition of $X_{\beta}$ (see Introduction) and
the fact that the positive root $\beta_{0}$ of the equation $x^{3}=x+1$ is the
smallest Pisot number \cite{DuPi}. Indeed, $\beta_{0}^{5}=\beta_{0}^{4}+1$ and
$X_{\beta_{0}}$ is a subshift of finite type, namely,
\[
X_{\beta_{0}}=\left\{  \bar{\varepsilon}\in\prod_{-\infty}^{+\infty
}\{0,1\}\mid\varepsilon_{n}=1\Rightarrow\varepsilon_{n+1}=\varepsilon
_{n+2}=\varepsilon_{n+3}=\varepsilon_{n+4}=0\right\}  .
\]
Now the desired claim follows from \cite[Lemma~3]{Pa} asserting that if
$\beta^{\prime}<\beta$, then $(d_{1}(\beta^{\prime}),d_{2}(\beta^{\prime
}),\dots)\prec(d_{1}(\beta),d_{2}(\beta),\dots)$.
\end{proof}

Let
\begin{align*}
\mathfrak{A}  &  =\{\bar{\varepsilon}\in X_{\beta}^{+}\mid\exists
n\in\mathbb{N}:\forall\alpha\in Z_{\beta},\;\bar{\alpha}+(\varepsilon
_{1},\ldots,\varepsilon_{n})\in Fin(\beta)\},\\
\mathfrak{A}_{n}  &  =\{\bar{\varepsilon}\in X_{\beta}^{+}\mid\forall\alpha\in
Z_{\beta},\;\bar{\alpha}+(\varepsilon_{1},\ldots,\varepsilon_{n})\in
Fin(\beta)\},\\
\mathfrak{A}^{\prime}  &  =\{\bar{\varepsilon}\in X_{\beta}^{+}\mid\exists
n\in\mathbb{N}:\varepsilon_{n+1}=\cdots=\varepsilon_{n+L+4}=0\}.
\end{align*}
We will write $tail(\bar{\varepsilon})=tail(\bar{\varepsilon}^{\prime})$ if
there exists $n\in\mathbb{N}$ such that $\varepsilon_{k}=\varepsilon
_{k}^{\prime},\;k\geq n$. The meaning of the above definitions consists in the
fact that if $\bar{\varepsilon}\in\mathfrak{A}\cap\mathfrak{A}^{\prime}$, then
$\bar{\varepsilon}\in\mathfrak{A}_{n}\cap\mathfrak{A}^{\prime}$ for some
$n\geq1$ and by the theorem from \cite{FrSo} mentioned above, $(\varepsilon
_{1},\ldots,\varepsilon_{n})+\bar{\alpha}=(\varepsilon_{1}^{\prime}%
,\ldots,\varepsilon_{n+L}^{\prime})$, whence by Lemma~\ref{glue}
\[
tail(\bar{\varepsilon}+\bar{\alpha})=tail(\bar{\varepsilon})
\]
(more precisely, the tail will stay unchanged starting with the $(n+L+1)$'th
symbol). It is obvious that $\mathfrak{A}=\cup_{n}\mathfrak{A}_{n}$. We wish
to prove that $\mu_{\beta}^{+}(\mathfrak{A}\cap\mathfrak{A}^{\prime})=1$. By
the ergodicity of $(X_{\beta}^{+},\mu_{\beta}^{+},\sigma_{\beta}^{+})$, we
have $\mu_{\beta}^{+}(\mathfrak{A}^{\prime})=1$, it suffices to show that
$\mu_{\beta}^{+}(\mathfrak{A})=1$. Let $\mathfrak{B}_{n}=X_{\beta}%
^{+}\setminus\mathfrak{A}_{n}$.

\begin{proposition}
\label{16}There exists a constant $\gamma=\gamma(\beta)\in(0,1)$ such that
\begin{equation}
\mu_{\beta}^{+}\left(  \bigcap_{k=1}^{n}\mathfrak{B}_{k}\right)  \leq
\gamma^{n}. \label{gamma}%
\end{equation}
\end{proposition}

\begin{proof}
We have
\begin{align*}
\mu_{\beta}^{+}(\mathfrak{B}_{1}\cap\mathfrak{B}_{2}\cap\ldots\cap\mathfrak
{B}_{n})  &  =\mu_{\beta}^{+}(\mathfrak{B}_{1})\cdot\prod_{k=2}^{n}\mu_{\beta
}^{+} (\mathfrak{B}_{k}\mid\mathfrak{B}_{k-1}\cap\ldots\cap\mathfrak{B}_{1})\\
&  \leq\prod_{k=2}^{n}\mu_{\beta}^{+}(\mathfrak{B}_{k}\mid\mathfrak{B}%
_{k-1}\cap\ldots\cap\mathfrak{B}_{1}).
\end{align*}
Since
\begin{align*}
\prod_{j=k-L}^{k}\mu_{\beta}^{+}(\mathfrak{B}_{j}\mid\mathfrak{B}_{j-1}%
\cap\ldots\cap\mathfrak{B}_{1})  &  =\frac{\mu_{\beta}^{+}(\mathfrak{B}%
_{k}\cap\ldots\cap\mathfrak{B}_{1})}{\mu_{\beta}^{+}(\mathfrak{B}_{k-L-1}%
\cap\ldots\cap\mathfrak{B}_{1})}\\
&  \leq\frac{\mu_{\beta}^{+}(\mathfrak{B}_{k}\cap\mathfrak{B}_{k-L-1}%
\cap\mathfrak{B}_{k-L-2}\cap\ldots\cap\mathfrak{B}_{1})}{\mu_{\beta}^{+}
(\mathfrak{B}_{k-L-1}\cap\ldots\cap\mathfrak{B}_{1})}\\
&  =\mu_{\beta}^{+}(\mathfrak{B}_{k}\mid\mathfrak{B}_{k-L-1}\cap\mathfrak
{B}_{k-L-2} \cap\ldots\cap\mathfrak{B}_{1}),
\end{align*}
we have
\begin{equation}
\mu_{\beta}^{+}(\mathfrak{B}_{1}\cap\mathfrak{B}_{2}\cap\ldots\cap\mathfrak
{B}_{n})\leq\prod_{k=2}^{[n/L]}\mu_{\beta}^{+}(\mathfrak{B}_{Lk}\mid
\mathfrak{B}_{Lk-L-1}\cap\mathfrak{B}_{Lk-L-2}\cap\ldots\cap\mathfrak{B}_{1}).
\label{n/L}%
\end{equation}
Now by the formula (\ref{C_1}), $\beta$ being weakly finitary (see (\ref{WF}))
and the deifnition of $L$ we have%
\[
\mu_{\beta}^{+}(\mathfrak{A}_{k+L}\mid\varepsilon_{k}=i_{k},\ldots
,\varepsilon_{1}=i_{1})\geq C_{2}>0
\]
for any $(i_{1},\ldots,i_{k})\in X_{\beta}^{(k)}$. Hence
\[
\mu_{\beta}^{+}(\mathfrak{B}_{Lk}\mid\mathfrak{B}_{Lk-L-1}\cap\mathfrak
{B}_{Lk-L-2}\cap\ldots\cap\mathfrak{B}_{1})\leq1-C_{2},
\]
and from (\ref{n/L}) we finally obtain the estimate%
\[
\mu_{\beta}^{+}(\mathfrak{B}_{1}\cap\mathfrak{B}_{2}\cap\ldots\cap\mathfrak
{B}_{n})\leq(1-C_{2})^{n/L},
\]
whence one can take $\gamma=(1-C_{2})^{1/L}$, and (\ref{gamma}) is proven.
\end{proof}

As a consequence we obtain the following claim about the irrational rotations
of the circle by the elements of $\mathbb{Z}[\beta]$. Let, as above,
$\bar{\alpha}$ denote the $\beta$-expansion of $\alpha$.

\begin{theorem}
For a weakly finitary Pisot unit $\beta$ and any $\alpha\in\mathbb{Z}%
[\beta]\cap\lbrack0,1)$ we have%
\[
tail(\bar{\varepsilon}+\bar{\alpha})=tail(\bar{\varepsilon})
\]
for $\mu_{\beta}^{+}$-a.e. $\bar{\varepsilon}\in X_{\beta}^{+}$.
\end{theorem}

\begin{proof}
We showed that $\mu_{\beta}^{+}(\cap_{1}^{\infty}\mathfrak{B}_{n})=0$, whence
$\mu_{\beta}^{+}(\mathfrak{A)=
\mu}_{\beta}^{+}(\cup_{1}^{\infty}\mathfrak{A}_{n})=1$.
\end{proof}

\noindent\textbf{Conclusion of the proof of Theorem~\ref{main}.} Fix
$k\in\mathbb{N}$. To complete the proof of Theorem~\ref{main}, it suffices to
show that the set
\[
\mathcal{D}^{(k)}=\{\bar{\varepsilon}\in X_{\beta}\mid\varepsilon_{j}%
\equiv\varepsilon_{j}^{\prime},\;j\geq k\;\forall\bar{\varepsilon}^{\prime
}\sim\bar{\varepsilon}\},
\]
has the full measure $\mu_{\beta}$. By Proposition \ref{16}, for
\[
\mathcal{D}_{N}^{(k)}=\{(\varepsilon_{-N},\varepsilon_{-N+1},\ldots)\in
X_{\beta}^{+}\mid(\ldots0,0,\varepsilon_{-N},\varepsilon_{-N+1},\ldots
)\in\mathcal{D}^{(k)}\},
\]
$\mu_{\beta}^{+}(\mathcal{D}_{N}^{(k)})\geq1-\gamma^{k-N}\rightarrow1$ as
$N\rightarrow+\infty$. Hence%
\[
\mu_{\beta}(\mathcal{D}^{(k)})=\lim_{N\rightarrow+\infty}\mu_{\beta}%
^{+}(\mathcal{D}_{N}^{(k)})=1,
\]
and therefore
\[
\mu_{\beta}\left(  \bigcap_{k=1}^{\infty}\mathcal{D}^{(k)}\right)  =1
\]
which implies (\ref{mu=1}). We have thus shown that for $\mu_{\beta}$-a.e.
$\bar{\varepsilon}\in X_{\beta}$, $\#\varphi_{\mathbf{t}}^{-1}(\varphi
_{\mathbf{t}}(\bar{\varepsilon}))=1$. Let $\mathcal{L}$ denote the image of
$\mu_{\beta}$ under $\varphi_{\mathbf{t}}$. Since $\mu_{\beta}$ is ergodic, so
is $\mathcal{L}$ and since $h_{\mu_{\beta}}(\sigma_{\beta})=\log\beta$, we
have $h_{\mathcal{L}}(T)=\log\beta$ as well. Hence $\mathcal{L=L}_{m}$ is the
Haar measure on the torus, as it is the unique ergodic measure of maximal
entropy. So, we proved that%
\[
\mathcal{L}_{m}\{x\in\mathbb{T}^{m}\mid\#\varphi_{\mathbf{t}}^{-1}%
(x)=1\}=1\text{,}%
\]
which is the claim of Theorem~\ref{main}.

As a corollary we obtain the following claim about the arithmetic
structure of $X_\beta$ itself.

\begin{proposition} Let $\sim$ denote the equivalence relation on $X_\beta$
defined by (\ref{er}) and $X_\beta':=X_\beta/\sim$. Then $X_\beta'$ is
a group isomorphic to $\mathbb{T}^m$.
\end{proposition}

Thus, $X_\beta$ is an {\em almost group} in the sense that
it suffices to ``glue" some $k$-tuples (for $k<\infty$)
within the set of measure zero to turn the two-sided $\beta$-compactum
for a weakly finitary Pisot unit $\beta$
into a group (which will be isomorphic to the torus of the corresponding
dimension). Note that in dimension~2 this factorisation can be
described more explicitly -- see \cite[Section~1]{SV98}.

The following claim is a generalisation of Theorem~4 from \cite{SV-Usp}. Let
$\mathcal{D}(T)$ denotes the \emph{centraliser }for $T$, i.e.,
\[
\mathcal{D}(T)=\{A\in GL(m,\mathbb{Z}):AT=TA\}.
\]

\begin{proposition}
For a Pisot automorphism whose matrix is algebraically conjugate to the
corresponding companion matrix there is a one-to-one correspondence between
the following sets:

\begin{enumerate}
\item  the fundamental homoclinic points for $T$;

\item  the bijective arithmetic codings for $T$;

\item  the units of the ring $\mathbb{Z}[\beta]$;

\item  the centraliser for $T$;
\end{enumerate}
\end{proposition}

\begin{proof}
We already know that any bijective arithmetic coding is parametrised by a
fundamental homoclinic point. Let $\mathbf{t}_{0}$ be such a point for $T$;
then any other fundamental homoclinic point $\mathbf{t}$ satsifies
$\mathbf{s}=u\mathbf{s}_{0}$, where $\mathbf{s}_{0}$ and $\mathbf{s}$ are the
corresponding $\mathbb{R}^{m}$-coordinates and $u\in\mathcal{U}_\beta$
-- the proof is essentially the same as in
Lemma~\ref{11}. On the other hand, if $\varphi_{\mathbf{t}}$ is a bijective
arithmetic coding for $T$, then as easy to see, $A:=\varphi_{\mathbf{t}%
}\varphi_{\mathbf{t}_{0}}^{-1}$ is a toral automorphism commuting with $T$
(this mapping is well defined almost everywhere on the torus, hence it can be
defined everywhere by continuity). Finally, if $u\in\mathcal{U}_{\beta}$ and
$u=\sum_{j=1}^{m-1}u_{j}\beta^{j}$, then $A:=\sum_{j=1}^{m-1}u_{j}M^{j}$
belongs to $GL(m,\mathbb{Z})$ and commutes with $M$ and vice versa.
\end{proof}

\medskip

\noindent\textbf{Example 1.} (see \cite{SV97}) Let $T$ be given by the matrix
$M=%
\begin{pmatrix}
1 & 1\\
1 & 0
\end{pmatrix}
$. Here $\beta$ is the golden ratio, and. $\xi_{0}=\frac{1}{\sqrt{5}}%
=\frac{-1+2\beta}{5}$, and
\[
\mathcal{U}_{\beta}=\{\pm\beta^{n},\ n\in\mathbb{Z}\}.
\]
Any bijective arithmetic coding for $T$ thus will be of the form
\[
\varphi(\overline{\varepsilon})=\lim_{N\rightarrow+\infty}\sum_{k=-N}%
^{+\infty}\varepsilon_{k}\beta^{-k}%
\begin{pmatrix}
\vartheta\beta^{n}/\sqrt{5}\\
\vartheta\beta^{n-1}/\sqrt{5}%
\end{pmatrix}
\operatorname{mod}\mathbb{Z}^{2},
\]
where $\vartheta\in\{\pm1\}$ and $n\in\mathbb{Z}$.

For more two-dimensional examples see \cite{SV98}.

\medskip

\noindent\textbf{Example 2.} Let $T$ be given by the matrix $M=%
\begin{pmatrix}
1 & 1 & 1\\
1 & 0 & 0\\
0 & 1 & 0
\end{pmatrix}
$. Here $\beta$ is the real root of the ``tribonacci'' equation
$x^3=x^2+x+1$; as is well known, $\beta$ is finitary in this case
(see, e.g., \cite{FrSo}).
We have $\xi_{0}=\frac{1}{-1-2\beta+3\beta^{2}}=\frac
{1+9\beta-4\beta^{2}}{22}$, and since $\mathbb{Z[\beta]}$ is the maximal order
in the field $\mathbb{Q}(\beta)$ and both conjugates of $\beta$ are complex,
again
\[
\mathcal{U}_{\beta}=\{\pm\beta^{n},\ n\in\mathbb{Z}\}
\]
(recall that by Dirichlet's Theorem, $\mathcal{U}_{\beta}$ must be
``one-dimensional", see, e.g., \cite{BorShaf}). Hence any bijective arithmetic
coding for $T$ is of the form
\[
\varphi(\overline{\varepsilon})=\lim_{N\rightarrow+\infty}\sum_{k=-N}%
^{+\infty}\varepsilon_{k}\beta^{-k}%
\begin{pmatrix}
\vartheta\frac{1+9\beta-4\beta^{2}}{22}\;\beta^{n}\\
\vartheta\frac{1+9\beta-4\beta^{2}}{22}\;\beta^{n-1}\\
\vartheta\;\frac{1+9\beta-4\beta^{2}}{22}\beta^{n-2}%
\end{pmatrix}
\operatorname{mod}\mathbb{Z}^{3},
\]
where $\vartheta\in\{\pm1\}$ and $n\in\mathbb{Z}$.

\medskip

\noindent\textbf{Example 3.} Let $M=%
\begin{pmatrix}
3 & 4 & 1\\
1 & 0 & 0\\
0 & 1 & 0
\end{pmatrix}
$. Here $\beta$ is the positive root of $x^{3}=3x^{2}+4x+1$. By the result
from \cite{Ak-cubic}, $\beta$ is finitary (see Introduction for the
definition) and the fundamental units of the ring are $\beta$ and
$3+\beta^{-1}$, i.e.,%
\[
\mathcal{U}_{\beta}=\{\pm\beta^{n},\pm(3+\beta^{-1})^{n},\ n\in\mathbb{Z}\}.
\]
Besides, $\xi_{0}=\frac{1}{3\beta^{2}-6\beta-4}=\frac{-13-21\beta+6\beta^{2}%
}{7}$. Hence any bijective arithmetic coding is either%
\[
\varphi(\overline{\varepsilon})=\lim_{N\rightarrow+\infty}\sum_{k=-N}%
^{+\infty}\varepsilon_{k}\beta^{-k}%
\begin{pmatrix}
\vartheta\frac{-13-21\beta+6\beta^{2}}{7}\;\beta^{n}\\
\vartheta\frac{-13-21\beta+6\beta^{2}}{7}\;\beta^{n-1}\\
\vartheta\;\frac{-13-21\beta+6\beta^{2}}{7}\beta^{n-2}%
\end{pmatrix}
\operatorname{mod}\mathbb{Z}^{3}%
\]
or%
\[
\varphi(\overline{\varepsilon})=\lim_{N\rightarrow+\infty}\sum_{k=-N}%
^{+\infty}\varepsilon_{k}\beta^{-k}%
\begin{pmatrix}
\vartheta\frac{-13-21\beta+6\beta^{2}}{7}(3+\beta^{-1})^{n}\\
\vartheta\frac{-13-21\beta+6\beta^{2}}{7}(3+\beta^{-1})^{n-1}\\
\vartheta\;\frac{-13-21\beta+6\beta^{2}}{7}(3+\beta^{-1})^{n-2}%
\end{pmatrix}
\operatorname{mod}\mathbb{Z}^{3},
\]
where $\vartheta\in\{\pm1\}$ and $n\in\mathbb{Z}$.

\medskip

\noindent\textbf{Example 4.} Finally, let $M=\left(
\begin{array}
[c]{cccc}%
1 & 0 & 0 & 1\\
1 & 0 & 0 & 0\\
0 & 1 & 0 & 0\\
0 & 0 & 1 & 0
\end{array}
\right) $. Here $\beta$ satisfies $x^{4}=x^{3}+1$. Let us show that $\beta$ is
weakly finitary. A direct inspection shows that the only nonzero tail for the
positive elements of $\mathbb{Z}[\beta]$ is $\overline{10000}$. Hence
$Z_{\beta}=\{0,\beta^{-2}+\beta^{-3},\beta^{-3}+\beta^{-4},\beta^{-4}%
+\beta^{-5},\beta^{-5}+\beta^{-6},\beta^{-6}+\beta^{-7}\}$. Let, for example,
$x=\beta^{-2}+\beta^{-3}$; since $x+\beta^{-5}=\beta^{-1}+\beta^{-3}%
=\beta^{-1}+\beta^{-4}+\beta^{-7}=1+\beta^{-7}\in Fin(\beta)$, we have by
periodicity $x+\beta^{-5n}\in Fin(\beta)$ for any $n\geq1$. The other cases of
$\alpha\in Z_{\beta}$ are similar. Hence $\beta$ is weakly finitary and we can
apply Theorem \ref{main}. It suffices to compute $\mathcal{U}_{\beta}$; by the
Dirichlet Theorem, it must be ``two-dimensional'' and it is easy to guess that
the second fundamental unit (besides $\beta$ itself) is $1+\beta$. Hence
$\mathcal{U}_{\beta}=\{\pm\beta^{n},\pm(1+\beta)^{n},\ n\in\mathbb{Z}\}$ and
the formula for a bijective arithmetic coding can be derived similarly to the
previous examples in view of $\xi_{0}=\frac{1}{-3\beta^{2}+4\beta^{3}}%
=\frac{-12-16\beta+73\beta^{2}+9\beta^{3}}{283}$.

\section{General arithmetic codings and related algebraic results}

In this section we will present some results for the case when $\mathbf{t}$ is
not necessarily fundamental or $T$ is not algebraically conjugate to the
companion matrix automorphism. We will still assume $\beta$ to be weakly
finitary. Let us begin with the case $T=T_{\beta}$ with a general $\mathbf{t}%
$. We recall that there is an isomorphism between the homoclinic group
$\mathcal{H}(T)$ and the group $\mathcal{P}_{\beta}$, i.e., $\mathbf{t}%
\leftrightarrow\xi$. Let $\varphi_{\xi}:X_{\beta}\rightarrow\mathbb{T}^{m}$ be
defined as usual:
\[
\varphi_{\xi}(\bar{\varepsilon})=\varphi_{\mathbf{t}}(\bar{\varepsilon}%
)=\sum_{k\in Z}\varepsilon_{k}T_{\beta}^{-k}\mathbf{t}=\lim_{N\rightarrow
+\infty}\left(  \sum_{k=-N}^{\infty}\varepsilon_{k}\beta^{-k}\right)
\begin{pmatrix}
\xi\\
\xi\beta^{-1}\\
\vdots\\
\xi\beta^{-m+1}%
\end{pmatrix}
\operatorname{mod}\mathbb{Z}^{m},
\]
where $\xi=\xi(\mathbf{t})\in\mathcal{P}_{\beta}$. The question is, what will
be the value of $\#\varphi_{\xi}^{-1}(x)$ for a $\mathcal{L}_{m}$-typical
$x\in\mathbb{T}^{m}$?

The next assertion answers this question; it is a generalization of the
corresponding result proven in \cite{SV98} for $m=2$ and for a finitary
$\beta$ in \cite{Sid}. Let $D=D(\beta)$ denote the discriminant of $\beta$ in
the field extension $\mathbb{Q}(\beta)/\mathbb{Q}$, i.e., the product
$\prod_{i\neq j}(\beta_{i}-\beta_{j})^{2}$, where $\{\beta_{1}=\beta,\beta
_{2},\ldots,\beta_{m}\}$ are the Galois conjugates of $\beta$.

\begin{theorem}
For an a.e. $x\in\mathbb{T}^{m}$ with respect to the Haar measure,
\[
\#\varphi_{\xi}^{-1}(x)\equiv|DN(\xi)|,
\]
where $N(\cdot)$ denotes the norm of an element of the extension
$\mathbb{Q}(\beta)/\mathbb{Q}$.\label{preimages}
\end{theorem}

\begin{proof}
Let $\varphi_{0}$ denote the bijective arithmetic coding for $T_{\beta}$
parametrised by $\xi_{0}$ and $\ell:=\xi/\xi_{0}\in\mathbb{Z}[\beta]$. If
$\ell=\sum_{i=0}^{m-1}c_{i}\beta^{i}$, then one can consider the mapping
$A_{\xi}:=\varphi_{\xi}\varphi_{0}^{-1}:\mathbb{T}^{m}\rightarrow
\mathbb{T}^{m}$; it will be well defined on the dense set and we may extend it
to the whole torus. By the linearity of both maps, $A_{\xi}$ is a toral
endomorphism. Thus, we have%
\begin{equation}
\varphi_{\xi}=A_{\xi}\varphi_{0}. \label{phi-xi}%
\end{equation}
Let $A^{\prime}_{\xi}$ is given by the formula $A^{\prime}_{\xi}=\sum
_{i=0}^{m-1}c_{i}T_{\beta}^{i}$. For the basis sequence $f^{(0)}%
=(\dots,0,0,\dots,0,1,0,\dots,0,0,\dots)$ with the unity at the first
coordinate we have%
\[
\begin{aligned}
(A_\xi\varphi_0)(f^{(0)})&=A_\xi(\xi_0,\xi_0\beta^{-1},\dots,\xi_0\beta
^{-m+1})
\operatorname{mod}\mathbb{Z}^m \\
&=\sum_{i=0}^{m-1}c_iT^i(\xi_0,\xi_0\beta^{-1},\dots,\xi_0\beta^{-m+1})
\operatorname{mod}\mathbb{Z}^m\\
&=\sum_{i=0}^{m-1}c_i\beta^i(\xi_0,\xi_0\beta^{-1},\dots,\xi_0\beta^{-m+1})
\operatorname{mod}\mathbb{Z}^m \\
&=(\xi,\xi\beta^{-1},\dots,\xi\beta^{-m+1})
\operatorname{mod}\mathbb{Z}^m=\varphi
_\xi(f^{(0)}).
\end{aligned}
\]
Therefore, by the linearity and continuty, we have $A_{\xi}=A^{\prime}_{\xi
}=\sum_{i=0}^{m-1}c_{i}T_{\beta}^{i}$. As $\varphi_{0}$ is 1-to-1 a.e.,
$\varphi_{\xi}$ will be $K$-to-1 a.e. with $K=|\det A_{\xi}|$. By definition,
$N(\ell)$ is the determinant of the matrix of the multiplication operator
$x\mapsto\ell x$ in the standard basis of $\mathbb{Q}(\beta)$, whence
$N(\ell)=\det A_{\xi}$, because $T_{\beta}$ is given by the companion matrix.
Finally, $N(\ell)=N(\xi)/N(\xi_{0})=DN(\xi)$, as by the result from
\cite[Section~2.7]{Sam}, $N(\xi_{0})=1/D$ whenever $\xi_{0}$ is as in
formula~(\ref{P_beta}).
\end{proof}

Note that historically the first attempt to find an arithmetic coding for a
Pisot automorphism was undertaken in \cite{Ber}. The author considered the
case $T=T_{\beta}$ and $\mathbf{t}$ given by the $\mathbb{R}^{m}$-coordinate
$\mathbf{s}=(1,\beta^{-1},\dots,\beta^{-m+1})$. From the above theorem follows

\begin{corollary}
The mapping
\[
\varphi_{\mathbf{t}}(\bar{\varepsilon})=\lim_{N\rightarrow+\infty}\sum
_{k=-N}^{+\infty}\varepsilon_{k}\beta^{-k}\left(
\begin{array}
[c]{c}%
1\\
\beta^{-1}\\
\vdots\\
\beta^{-m+1}%
\end{array}
\right)  \operatorname{mod}\;\mathbb{Z}^{m}%
\]
is $|D|$-to-1 a.e.
\end{corollary}

Suppose now $T$ is not necessarily algebraically conjugate to $T_{\beta}$. Let
$M$ be, as usual, the matrix of $T$, and for $\mathbf{n}\in\mathbb{Z}^{m}$ the
matrix $B_{M}(\mathbf{n})$ be defined as follows (we write it column-wise):
\[
\begin{aligned}
B_{M}(\mathbf{n})=(&M\mathbf{n},(M^{2}-k_{1}M)\mathbf{n},
(M^{3}-k_{1}M^{2}-k_{2}M)\mathbf{n},\ldots,\\&M^{m-1}-k_{1}M^{m-2}-
\cdots-k_{m-2}M)\mathbf{n}, k_{m}\mathbf{n}).
\end{aligned}
\]

\begin{lemma}
\label{Semic}Any integral square matrix satisfying the relation%
\begin{equation}
BM_{\beta}=MB \label{semiconj}%
\end{equation}
is $B=B_{M}(\mathbf{n)}$ for some $\mathbf{n}\in\mathbb{Z}^{m}$.
\end{lemma}

\begin{proof}
Let $B$ be written column-wise as follows: $B=(\mathbf{n}_{1},\ldots
,\mathbf{n}_{m})$. Then by (\ref{semiconj}) and the definition of $M_{\beta}%
$,
\[
(k_{1}\mathbf{n}_{1}+\mathbf{n}_{2},k_{2}\mathbf{n}_{1}+\mathbf{n}_{3}%
,\ldots,k_{m-1}\mathbf{n}_{1}+\mathbf{n}_{m},k_{m}\mathbf{n}_{1}%
)=(M\mathbf{n}_{1},\ldots,M\mathbf{n}_{m}),
\]
whence by the fact that $k_{m}=\pm1$, we have $B=B_{M}(\mathbf{n})$ for
$\mathbf{n}=k_{m}\mathbf{n}_{m}$.
\end{proof}

\begin{definition}
The integral $m$-form of $m$ variables defined by the formula%
\[
f_{M}(\mathbf{n})=\det B_{M}(\mathbf{n})
\]
will be called the \textbf{form associated with} $T$.
\end{definition}

\begin{proposition}
\label{numpreimage} Let $\mathbf{t}\in\mathcal{H}(T)$. Then there exists
$\mathbf{n}\in$ $\mathbb{Z}^{m}$ such that%
\[
\#\varphi_{\mathbf{t}}^{-1}(x)\equiv|f_{M}(\mathbf{n})|
\]
for $\mathcal{L}_{m}$-a.e. point $x\in\mathbb{T}^{m}$. Hence the minimum of
the number of preimages for an arithmetic coding of a given automorphism $T$
equals the arithmetic minimum of the associated form $f_{M}$.
\end{proposition}

\begin{proof}
Let $\tilde{B}:=\varphi_{\mathbf{t}}\varphi_{0}^{-1}$, where $\varphi_{0}$ is
a certain bijective arithmetic coding for $T_{\beta}$. Then $\tilde{B}$ is a
linear mapping from $\mathbb{T}^{m}$ onto itself defined a.e.; let the same
letter denote the corresponding toral endomorphism. Then $\tilde{B}T_{\beta
}=\varphi_{\mathbf{t}}\varphi_{0}^{-1}T_{\beta}=\varphi_{\mathbf{t}}%
\sigma_{\beta}\varphi_{0}^{-1}=T\varphi_{\mathbf{t}}\varphi_{0}^{-1}%
=T\tilde{B}$. Therefore the matrix $B$ of the endomorphism $\tilde{B}$
satisfies (\ref{semiconj}), whence by Lemma \ref{Semic},
$B=B_{M}(\mathbf{n})$ for some $\mathbf{n}\in\mathbb{Z}^{m}$. Hence
$\varphi_{\mathbf{t}}=B_{M}(\mathbf{n})\varphi_{0}$, and we are done.
\end{proof}

\begin{remark}
It would be helpful to know whether there is any relationship between the
$\mathbf{n}$ in the proposition and the $\mathbb{Z}^{m}$-coordinate of
$\mathbf{t}$.
\end{remark}

\begin{theorem}
\label{bac} The following conditions are equivalent:

\begin{enumerate}
\item  An automorphism $T$ admits\ a bijective arithmetic coding.

\item $T$ is algebraically conjugate to $T_{\beta}$.

\item  The equation
\[
f_{M}(\mathbf{n})=\pm1
\]
has a solution in $\mathbf{n}\in\mathbb{Z}^{m}$.

\item  There exists a homoclinic point $\mathbf{t}$ such that for its
$\mathbb{Z}^{m}$-coordinate $\mathbf{n}$,%
\[
\left\langle M^{k}\mathbf{n}\mid k\in\mathbb{Z}\right\rangle =\mathbb{Z}^{m}.
\]
\end{enumerate}
\end{theorem}

\begin{proof}
(2)$\Rightarrow$(1): see Remark \ref{conj};

(1)$\Rightarrow$(2): see the Proposition~\ref{numpreimage};

(2)$\Leftrightarrow$(3): also follows from Proposition~\ref{numpreimage};

(2)$\Leftrightarrow$(4): it is obvious that $M_{\beta}$ satisfies this
property (take $\mathbf{n}=(0,0,\ldots,0,1)$). Hence so does any $M$ which is
conjugate to $M_{\beta}$.
\end{proof}

Recall that a matrix $M\in GL(m,\mathbb{Z})$ is called \textit{primitive} if
there is no matrix $K\in GL(m,\mathbb{Z})$ such that $M=K^{n}$ for $n\geq2$.
Following \cite{SV98}, we ask the following question: can a Pisot toral
automorphism given by a non-primitive matrix admit a bijective arithmetic
coding?

Note first that one can simplify the formula for $f_{M}$. Namely, since the
determinant of a matrix stays unchanged if we multiply one column by some
number and add to another column, we have
\begin{equation}
f_{M}(\mathbf{n})=\pm\det\;(\mathbf{n},M\mathbf{n},\ldots,M^{m-1}%
\mathbf{n}).\label{fM}%
\end{equation}

\begin{proposition}
There exists a sequence of integers $\mathcal{N}_{n}(\beta)$ such that%
\[
f_{M^{n}}=\pm\mathcal{N}_{n}(\beta)\cdot f_{M}.
\]
More precisely,
\[
\mathcal{N}_{n}(\beta)=\det\left(
\begin{array}
[c]{ccc}%
a_{n}^{(1)} & \ldots &  a_{n}^{(m)}\\
a_{2n}^{(1)} & \ldots &  a_{2n}^{(m)}\\
\vdots & \ddots & \vdots\\
a_{(m-1)n}^{(1)} & \ldots &  a_{(m-1)n}^{(m)}%
\end{array}
\right)  ,
\]
where $\{a_{n}^{(j)}\}_{j=1}^{m}$ are defined as the coefficients of the
equation%
\[
\beta^{n}=a_{n}^{(1)}\beta^{m-1}+a_{n}^{(2)}\beta^{m-2}+\cdots+a_{n}%
^{(m-1)}\beta+a_{n}^{(m)}.
\]
\end{proposition}

\begin{proof}
By (\ref{fM}), the definition of $a_{n}^{(j)}$ and the Hamilton-Cayley
Theorem,%
\begin{align*}
f_{M^{n}}(\mathbf{n})  & =\pm\det\;(\mathbf{n},M^{n}\mathbf{n},M^{2n}%
\mathbf{n},\ldots,M^{(m-1)n}\mathbf{n})\\
& =\pm\det\;\left(  \mathbf{n},\left(  \sum_{j=1}^{m}a_{n}^{(j)}%
M^{m-j}\right)  \mathbf{n},
\ldots,\left(  \sum_{j=1}^{m}a_{(m-1)n}^{(j)}M^{m-j}\right)
\mathbf{n}\right)  =\\
& =\pm\mathcal{N}_{n}(\beta)\cdot\det\;(\mathbf{n},M\mathbf{n},M^{2}%
\mathbf{n},\ldots,M^{n}\mathbf{n}).
\end{align*}
\end{proof}

\begin{corollary}
A non-primitive matrix $M^{n}\in GL(m,Z)$ is algebraically conjugate to the
corresponding companion matrix if and only if so is $M$,
and $\mathcal{N}_{n}(\beta)=\pm1$.
\end{corollary}

Let us deduce some corollaries for smaller dimensions.

\begin{corollary}
(see \cite{SV98}) For $m=2$ the automorphism given by a non-primitive matrix
$M^{n}$ admits a bijective arithmetic coding if and only if $n=2$ and
$Tr(M)=\pm1$.
\end{corollary}

\begin{corollary}
For $m=3$ the matrix $K=M^{2},M\in GL(3,\mathbb{Z})$, is algebraically
conjugate to the corresponding companion matrix if and only if $\beta$
satisfies one of the following equations:

\begin{enumerate}
\item $\beta^{3}=r\beta^{2}+1,\;r\geq1$;

\item $\beta^{3}=r\beta^{2}-1,\;r\geq3$;

\item $\beta^{3}=2\beta^{2}-\beta+1$.
\end{enumerate}
\end{corollary}

\begin{proof}
We have $\mathcal{N}_{2}(\beta)=\det\left(
\begin{array}
[c]{cc}%
1 & k_{1}^{2}+k_{2}\\
0 & k_{1}k_{2}+k_{3}%
\end{array}
\right)  =k_{1}k_{2}+k_{3}=\pm1$. The case $k_{3}=+1$ thus leads to subcases 1
and 3 and $k_{3}=-1$ yields subcase 2.
\end{proof}

Note that if $M$ is the matrix for the ``tribonacci automorphism'' (see
Example~2), then apparently the only power of $M$ that is algebraically
conjugate to the corresponding companion matrix, is the cube! Indeed,
$\mathcal{N}_{2}(\beta)=2,\mathcal{N}_{3}(\beta)=-1,\mathcal{N}_{4}%
(\beta)=-8,\mathcal{N}_{5}(\beta)=29,$ etc. It seems to be an easy exercise to
prove this rigoriously; we leave it to the reader.

\noindent\textbf{Example 5.} Let $M=%
\begin{pmatrix}
1 & 1 & 0\\
2 & 3 & 1\\
1 & 1 & 1
\end{pmatrix}
$. Here $\beta$ satisfies $x^{3}=5x^{2}-4x+1$ and the form associated with $M$
is (we write $\mathbf{n}=(x,y,z)^{\prime}$)%
\[
f_{M}(x,y,z)=x^{3}+2x^{2}z-xy^{2}-xyz+3xz^{2}+y^{3}-3y^{2}z+2yz^{2}+z^{3}.
\]
Obviously, the Diophantine equation $f_{M}(x,y,z)=\pm1$ has a solution,
namely, $x=1,y=z=0$. Hence by Theorem~\ref{bac}, $M$ is algebraically
conjugate to $M_{\beta}$; for example, $B=B_{M}(1,0,0)=%
\begin{pmatrix}
1 & 2 & -1\\
2 & -1 & 0\\
1 & -1 & 0
\end{pmatrix}
$ conjugates them. To show that $T$ admits a bijective arithmetic coding, it
suffices to check that $\beta$ is weakly finitary. We leave it to the
interested reader.

In \cite{SV98} the author together with A.~Vershik considered the case $m=2$.
Here if $M=\left(
\begin{array}
[c]{cc}%
a & b\\
c & d
\end{array}
\right)  $, then for $\sigma=\det M=\pm1$,%
\[
f_{M}(x,y)=\left|
\begin{array}
[c]{cc}%
ax+by & -\sigma x\\
cx+dy & -\sigma y
\end{array}
\right|  =\sigma(cx^{2}-(a-d)xy-by^{2}),
\]
and we related the problem of arithmetic codings to the classical theory of
binary quadratic forms. In particular, $T$ admits a bijective arithmetic
coding if and only if the Diophantine equation
\[
cx^{2}-(a-d)xy-by^{2}=\pm1
\]
is solvable.

The theory of general $m$-forms of $m$ variables does not seem to be well
developed; nonetheless, we would like to mention a certain algebraic result
which looks relevant. Recall that two integral forms are called
\textit{equivalent} if there exists a unimodular integral change of variables
turning one form into another.

\begin{proposition}
Let $M_{1},M_{2}$ in $GL(m,\mathbb{Z)}$ be conjugate, and%
\[
AM_{1}A^{-1}=M_{2},
\]
where $A\in GL(m,\mathbb{Z)}$. Then $f_{M_{1}}$ is equivalent either to
$f_{M_{2}}$ or to $-$ $f_{M_{2}}$, and moreover,%
\begin{equation}
A^{\prime}f_{M_{2}}A=\det A\cdot f_{M_{1}}, \label{detA}%
\end{equation}
where $A^{\prime}$ is the transpose of $A$ (we identify a form with the
symmetric matrix which defines it).
\end{proposition}

\begin{proof}
Since $M_{1}$ and $M_{2}$ are conjugate, they have one and the same
characteristic polynomial. By the definition of $f_{M}$ we have%
\begin{align*}
f_{M_{2}}(A\mathbf{v})  &  =\det(M_{2}A\mathbf{v},(M_{2}^{2}-k_{1}%
M_{2})A\mathbf{v},\ldots,A\mathbf{v})\\
&  =\det(AM_{1}\mathbf{v},A(M_{1}^{2}-k_{1}M_{1})\mathbf{v},\ldots
,A\mathbf{v})\\
&  =\det A\cdot f_{M_{1}}(\mathbf{v}),
\end{align*}
which is equivalent to (\ref{detA}).
\end{proof}

\end{document}